\documentclass{amsart}
\usepackage{amsmath,amsthm,amssymb,IMjournal}

\usepackage{fancybox,float,graphicx,latexsym,subfigure,srcltx,times,tikz,url,amsthm}

\newcommand{\C}{\mathbb{C}}
\newcommand{\R}{\mathbb{R}}
\newcommand{\Q}{\mathbb{Q}}
\newcommand{\N}{\mathbb{N}}
\newcommand{\Z}{\mathbb{Z}}

\newcommand{\BB}{\mathcal{B}}
\newcommand{\CC}{\mathcal{C}}
\newcommand{\FF}{\mathcal{F}}
\newcommand{\NN}{\mathcal{N}}
\newcommand{\MM}{\mathcal{M}}
\newcommand{\II}{\mathcal{I}}
\newcommand{\PP}{\mathcal{P}}
\renewcommand{\SS}{\mathcal{S}}

\newcommand{\col}{\operatorname{col}}
\newcommand{\dd}{\operatorname{d}}
\newcommand{\diag}{\operatorname{diag}}
\newcommand{\p}{\operatorname{p}}
\newcommand{\per}{\operatorname{per}}
\newcommand{\rk}{\operatorname{rk}}
\newcommand{\row}{\operatorname{row}}
\newcommand{\tr}{\operatorname{tr}}
\newcommand{\Trop}{\operatorname{Trop}}
\newcommand{\vv}{\operatorname{v}}

\begin{document}
\newtheorem{The}{Theorem}[section]

\newtheorem{lem}[The]{Lemma}
\newtheorem{cor}[The]{Corollary}
\newtheorem{prop}[The]{Proposition}
\newtheorem{ex}[The]{Example}
\newtheorem{rem}[The]{Remark}
\newtheorem{qst}[The]{Question}
\newtheorem{dfn}[The]{Definition}
\newtheorem{nota}[The]{Notation}
\newtheorem{con}[The]{Conjecture}

\numberwithin{equation}{section}

\title{Isocanted alcoved polytopes}

\author{\|Mar\'{\i}a Jes\'{u}s |de la Puente|, Madrid,
        \|Pedro Luis |Claver\'{\i}a|, Zaragoza}

%

\abstract
Through tropical normal idempotent matrices, we introduce isocanted alcoved polytopes, computing their $f$--vectors
and checking the
validity of  the following five conjectures: B\'{a}r\'{a}ny, unimodality, $3^d$, flag and cubical lower bound  (CLBC).
Isocanted alcoved polytopes are centrally symmetric, almost simple cubical polytopes. They are zonotopes. We show that, for
each dimension,  there is a unique combinatorial type. In dimension $d$, an isocanted alcoved polytope  has $2^{d+1}-2$ vertices, its
face lattice  is the lattice of proper subsets of $[d+1]$ and its diameter is $d+1$. They are realizations of $d$--elementary
cubical polytopes. The $f$--vector of a $d$--dimensional  isocanted alcoved polytope attains its maximum at the integer  $\lfloor d/3\rfloor$.
\endabstract

\keywords
cubical polytope, isocanted, alcoved, centrally symmetric, almost simple, zonotope,   $f$--vector, cubical $g$--vector, unimodal, flag, face lattice, log--concave sequence, tropical normal idempotent matrix, symmetric matrix.
\endkeywords

\subjclass
52B12, 15A80
\endsubjclass

\thanks
The first author is partially supported by  Ministerio de Econom\'{\i}a y Competitividad, Proyecto I+D MTM2016-76808-P, Ministerio de
Ciencia e Innovaci\'{o}n, Proyecto PID--2019--10770 GB-I00
and by UCM research group 910444.
\endthanks

\section{Introduction}\label{sec1}

This paper deals with  $f$--vectors of isocanted alcoved polytopes. A \emph{polytope} is the convex hull of a finite set of points
in $\R^d$. A polytope is a \emph{box} if its facets are only of one sort:  $x_i=cnst$, $i\in [d]$.
A polytope is \emph{alcoved} if its facets are only of two sorts:  $x_i=cnst$ and
$x_i-x_j=cnst$, $i,j\in[d]$, $i\neq j$.    Every alcoved polytope can be viewed as
the perturbation of  a box. In a box we distinguish two opposite vertices and the perturbation consists on canting
(i.e., beveling, meaning producing a flat face upon something) \label{beveling}
some (perhaps  all) of  the  $(d-2)$--faces of the box \emph{not meeting} the distinguished vertices. When the perturbation  happens for
all  such $(d-2)$--faces and with the same positive
magnitude, we obtain as a  result  an \emph{isocanted alcoved polytope}. The notion makes sense only  for $d\ge2$.

The $f$--vector of a $d$--polytope $\PP$ is the tuple  $(f_0,f_1,\ldots, f_{d-1})$, where $f_j$ is the number of  $j$--dimensional faces in $\PP$,  for $j=0,1,2,\ldots,{d-1}$. The $f$--vector can be extended with $f_{d}=1$.
It is well known that the $f$--vector of a $d$--box is
\begin{equation}\label{eqn:f_vector_box}
B_{d,j}=2^{d-j}{d\choose j},\quad j=0,1,\ldots,d.
\end{equation}
The quest for $f$--vectors is unrelenting.
As Ziegler writes in \cite{Ziegler} \lq\lq on some fundamental problems
embarrassingly little progress was made; one notable such problem concerns the
shapes of $f$--vectors\rq\rq and  \lq\lq new polytopes with interesting $f$--vectors should be produced\rq\rq and also
\lq\lq it seems that overall, we are short of examples.\rq\rq

The main result in this  paper is that  the $f$--vector of an isocanted $d$--alcoved polytope is given by
\begin{equation}\label{eqn:Idj}
I_{d,j}=(2^{d+1-j}-2){{d+1}\choose{j}},\ j=0,1,\ldots,d-1,\quad I_{d,d}=1.
\end{equation}
The numbers $I_{d,j}$ are even, for $j\le d-1$, because isocanted alcoved $d$--polytopes are centrally symmetric.
We verify several conjectures
for $f$--vectors, namely,  unimodality, B\'{a}r\'{a}ny, Kalai $3^d$ and flag conjectures as well as  CLCB.
Further properties are proved, showing that isocanted alcoved polytopes are $d$--elementary cubical, almost simple and zonotopes.

\medskip
The paper is organized as follows. In section  \ref{sec:isocanted} we give the definition and then, in Theorem \ref{thm:charact}, we prove
a crucial characterization: isocanted  alcoved polytopes are those
alcoved  polytopes having a unique vertex for each proper subset of $[d+1]$.  Concrete examples are given in Example \ref{ex:proof_thm}.
It follows from Theorem \ref{thm:charact} that the face lattice of an isocanted
alcoved $d$--polytope is the lattice of proper subsets of $[d+1]$. It is  proved that isocanted alcoved polytopes
are cubical and are zonotopes.
In section \ref{sec:d=3,4} we explain in detail the cases of dimensions 3 and 4, providing figures which help the reader visualize the
many properties of these polytopes. We compute two invariants of  4--isocanted alcoved polytopes: fatness and $f_{03}$.
In section \ref{sec:conjectures} we prove that the five mentioned conjectures hold true for isocanted alcoved polytopes.
Log--concavity provides a short proof of the   unimodality of $I_{d,j}$, for fixed $d\ge2$. We also prove that
the maximum of $I_{d,j}$ is attained at the integer
 $\lfloor\frac{d}{3}\rfloor$. We show that the diameter is $d+1$.

\medskip
This paper encompasses tropical matrices and classical polytopes, in the sense that tropical matrices are the means  to
describe certain polytopes.
We use several sorts of \emph{special matrices}, operated with tropical addition $\oplus=\max$  and tropical multiplication
$\odot=+$, such
as: normal idempotent (with respect to  $\odot$), visualized normal idempotent matrices, symmetric normal idempotent matrices and,  among these,
box matrices, cube matrices and isocanted matrices.

Tropical linear algebra and tropical algebraic geometry are fascinating, new,  fast growing areas
of mathematics with new and important results. For our purposes we recommend
\cite{Brugalle_fran,Brugalle_engl,Butkovic_Libro, Develin_Sturm,Litvinov_ed,Litvinov_ed_2,Mikhalkin_W,Richter,Speyer}
among many others. Alcoved polytopes have been first studied
in \cite{Lam_Postnikov,Werner_Yu}, then  in  \cite{Puente_kleene, Puente_QE}. Cubical polytopes have been addressed
in \cite{Adin,Adin_et_al,Blind_fewer,Blind_almost_simple,Jockusch}. General
references for polytopes
are \cite{Barvinok,Polytopes_Book,Grunbaum,Kalai_Ziegler_eds,Senechal_Libro,Schmitt_Ziegler,Ziegler_Libro,Ziegler}.
Normal idempotent matrices have been used in \cite{Puente_QE,Yu_Zhao_Zeng}. Idempotent matrices, also called Kleene stars,
have been used in \cite{Puente_kleene,Sergeev,Tran} in connection to polytopes.

\section{Background and notations}
\medskip
Well--known definitions and facts are presented here.
The set $\{1,2,\ldots,d+1\}$ is denoted $[d+1]$ and ${{[d+1]}\choose j}$ denotes the family of subsets of $[d+1]$ of cardinality
$j$. The origin in $\R^d$ is denoted  ${\bf{0}}$.  Maximum and minimum are taken componentwise in $\R^d$. A \emph{polyhedron}
in $\R^d$ is the intersection of a finite number of halfspaces.
It may be unbounded. A bounded polyhedron is called a  \emph{polytope} and every  polytope  is the convex hull of a finite set of points.
A $d$--polyhedron is a polyhedron of dimension $d$.
A  $d$--polyhedron $\PP$ is   \emph{alcoved} if its facets are only of two types: $x_i=cnst$ and
$x_i-x_j=cnst$, $i,j\in[d]$, $i\neq j$. A double index notation is useful  here because, in this way, we  can gather the
coefficients in a matrix over $\R\cup\{\pm\infty\}$: indeed,
write
\begin{equation}\label{eqn:ineq_ij}
a_{i,j}\le x_i-x_j\le -a_{j,i}
\end{equation} and, similarly,
\begin{equation}\label{eqn:ineq_i}
a_{i,d+1}\le x_i\le -a_{d+1,i}.
\end{equation} Then,
setting $a_{i,j}=\pm\infty$ if one facet $x_i-x_j=cnst$ is not specified, and  letting (by convention) $a_{i,i}=0$, for all
$i\in[d+1]$, we get a square matrix
$A=[a_{i,j}]\in M_{d+1}(\R\cup\{\pm\infty\})$ from $\PP$.  We write $\PP=\p(A)$ to express the former relationship between the
 polyhedron $\PP$ and the matrix $A$.
In addition to $a_{i,i}=0, i\in [d+1]$,\label{obs:null diagonal} the entries of the matrix $A$ satisfy
$-\infty\le a_{i,j}\le -a_{j,i}\le +\infty$, for all $i,j\in[d+1]$.
Different matrices $A$ may give rise to the same polyhedron.

\begin{dfn}[Alcoved polytope (AP)]\label{dfn:alcoved}
A $d$--polytope $\PP\subset \R^d$  is \emph{alcoved} if there exist
constants  $a_{i,j}\in \R$  such that $x\in\PP$ if and only if  $a_{i,d+1}\le x_i\le -a_{d+1,i}$, for all $i\in[d]$, and
$a_{i,j}\le x_i-x_j\le -a_{j,i}$, for all $i,j\in[d+1]$. Letting  $A=[a_{i,j}]\in M_{d+1}(\R)$,  we write $\PP=\p(A)$.
\end{dfn}

Important particular cases provide special matrices as follows:
\begin{enumerate}
\item ${\bf{0}}\in \p(A)$ if and only if  $A$ is \emph{normal} (N) (meaning  $a_{i,i}=0, a_{i,j}\le0, \forall i,j$), (see
\cite{Butkovic_Libro,Yu_Zhao_Zeng})
\item if ${\bf{0}}\in \p(A)$, then  $A$ describes $\p(A)$ \emph{optimally} (or tightly) if and only if $A$ is \emph{normal idempotent} (NI) (meaning that, in addition to normality, we have $A\odot A=A$, which requires that  $a_{i,j}+a_{jk}\le a_{ik}, \forall i,j,k$)\footnote{The family of normal idempotent matrices is a  subclass of the family of Kleene star matrices.} (see \cite{Puente_kleene,Sergeev,Tran}).
\item  for each alcoved polytope $\PP$ containing ${\bf{0}}$ there exist a unique NI
matrix $A$ such that $\PP=\p(A)$  (see Lemma 2.6 in \cite{Puente_kleene} and \cite{Sergeev,Tran}).  \label{item:2}
\end{enumerate}

Combinatorial properties of polytopes are, by nature,   translation invariant.
Every translate of an alcoved polytope is alcoved.  For each general alcoved polytope $\PP$,  infinitely many  translates $\PP'$ of
$\PP$  exist such that  ${\bf{0}}\in\PP'$. We can choose any such $\PP'$ to study $\PP$, and we know that $\PP'=\p(A)$ for a unique NI matrix
$A$. Often, we choose  $\PP'$ in two special locations with respect to  ${\bf{0}}$, each location  corresponding to a subclass of  NI matrices:
\begin{enumerate}
\item ${\bf{0}}=\max\p(A)$ if and only if $A$ is \emph{visualized normal idempotent} (VNI)  (in addition to NI,  the entries of $A$ satisfy $a_{d+1,i}=0,\forall i$),\label{item:a} (see
\cite{Butkovic_Libro,Puente_kleene,Puente_QE})
\item $\p(A)=-\p(A)$ if and only if $A$ is \emph{symmetric normal idempotent} (SNI) (in addition to NI, the entries of $A$ satisfy  $a_{i,j}=a_{j,i}, \forall i,j$), \label{item:b} (see
\cite{Jimenez_Puente,Puente_QE}).
\end{enumerate}
From \cite{Puente_QE}, we know that translation of an alcoved polyhedron $\p(A)$ corresponds to conjugation of its matrix $A$ by a
diagonal matrix (with null last diagonal entry). \label{comment:trans}

\medskip
Our aim is, after defining isocanted alcoved polytopes,  to compute the $f$--vector of those.
But, what is  already known about vertices of  an alcoved polytope $\p(A)$ in $\R^d$? First, the number of  vertices of $\p(A)$ is
bounded above by
${2d}\choose{d}$ and this bound is sharp (see \cite{Develin_Sturm,Speyer}). \emph{Which points are  vertices of} $\p(A)$?\label{qst:which}
In order to answer this question we introduce (a) the \emph{auxiliary matrix} $A_0$ and (b) the notion of \emph{tropical linear subspace} (by linear, we mean affine linear.)

For $A=[a_{i,j}]$,
the matrix
$A_0=[\alpha_{i,j}]$\label{dfn:A0} is defined by
$\alpha_{i,j}:=a_{i,j}-a_{d+1,j}=a_{i,j}\odot(-a_{d+1,j})$.\footnote{Notice that
$A_0$ might  be not normal.} The columns in $A_0$ are scalar multiples (with respect to  $\odot$) of the  columns in  $A$. The fact that $\diag(A)$
is zero  implies that $\row(d+1,A_0)$
is zero (and conversely), so that the columns in $A_0$ belong to the hyperplane  $\{x\in \R^{d+1}: x_{d+1}=0\}$
which is identified with $\R^d$.\footnote{This way of going  from $\R^{d+1}$ to $\R^d$, viewed as a hyperplane,  is analogous to going from projective to affine space, by intersecting with  the hyperplane $x_{d+1}=1$, in classical geometry.}
Besides, if $A$ is NI, then $A=A_0$ if and only if $A$ is VNI.

Inequalities  (\ref{eqn:ineq_i}) are transformed into
\begin{equation}
\alpha_{i,d+1}\le x_i\le \alpha_{i,i}, \ \forall i\in [d]
\end{equation} which yield the
following facts
\begin{equation}\label{eqn:min_and_max_PA}
\min\p(A)=\col(d+1,A_0),\qquad\max\p(A)=\diag(A_0).
\end{equation}
Besides, $\p(A)$ is the family of all tropical affine combinations of columns of $A_0$
(see Theorem 2.1\cite{Puente_kleene}, Proposition 12 \cite{Tran})\footnote{Here tropical geometry does not mimic classical geometry, since affine combinations do not produce the whole  tropical linear subspace, but only a bounded subset of it.}
\begin{equation}
\p(A)=\{x\in \R^{d+1}:x_{d+1}=0, x=\lambda_1\odot\col(1,A_0)\oplus\cdots\oplus\lambda_{d+1}\odot\col(d+1,A_0),
\lambda_j\in \R, 0=\lambda_1\oplus\cdots\oplus\lambda_{d+1}\}.
\end{equation}
$\p(A)$ is a proper subset of  the unique linear subspace determined by the columns of $A_0$.
In particular, the columns of $A_0$
are some of the vertices of $\p(A)$. They are  called the  \emph{generators} of  $\p(A)$. The rest of vertices of
$\p(A)$ are tropical linear combinations of the generators, and are thus called \emph{generated vertices} of  $\p(A)$.
In order to explain this, we must first define tropical linear subspaces.
A \emph{tropical linear subspace} is the \emph{tropicalization} of a linear subspace of $K^d$, where $K:=\C\{\{t\}\}$ is the field of
\emph{Puiseux series}. If $L\subseteq K^d$ is a linear subspace and
$I(L)\subseteq K[x_1^{\pm1},x_2^{\pm1},\ldots,x_d^{\pm1}]$ is the ideal of
all \emph{Laurent polynomials} vanishing on $L$, consider $q\in I(L)$,  $q=\sum_{s\in S}a_s\textbf{x}^s$, with
$\textbf{x}=(x_1,x_2,\ldots,x_d)$ variables, $s=(s_1,s_2,\ldots,s_d)\in S$ exponents,  $S\subset \N^d$ a finite set,
$a_s\in K$ and $\textbf{x}^s=x_1^{s_1}x_2^{s_2}\cdots x_d^{s_d}$.
Then, consider the tropicalization of $q$
\begin{equation}
\Trop(q):=\bigoplus_{s\in S}-\vv(a_s)\odot x_1{s_1}\odot x_2{s_2}\odot \cdots \odot x_d{s_d}=
\max_{s\in S}-\vv(a_s)+ x_1{s_1}+ x_2{s_2}+ \cdots +x_d{s_d}
\end{equation}
where tropical powers are transformed into products,  $\vv:K\setminus\{0\} \rightarrow \Q$ is the standard valuation
(i.e., the order of vanishing
of a series). The  \emph{corner locus} of  $\Trop(q)$ is, by definition,  the collection of
points  $x\in\R^d$  where the maximum in
$\Trop(q)(x)$ is attained, at
least, twice.\footnote{The translation to tropical mathematics of the
expression \lq\lq equal to zero \rq\rq or \lq\lq zero set \rq\rq is \lq\lq the maximum is attained, at least, twice.\rq\rq} Finally  $\Trop(L)$ is, by definition, the closure of the intersection of corner loci, for all $q\in I(L)$.
Since the \emph{corner locus} of  $\Trop(q)$  is piecewise linear, then tropical linear subspaces are polyhedral complexes.
\footnote{Unlike classical geometry, it is not true that, in $d$--dimensional space, the intersection of a  generic family of $(d-k)$ tropical linear hyperplanes is a tropical linear subspace of
dimension $k$.}

Notice that a unique tropical linear subspace is determined by
each  subset of generators (i.e., of columns of $A_0$). A convenient notation is $L_A(W)$, for each
$W\in {{[d+1]}\choose{j}}$ with $1\le j\le d$. We will write $L(W)$, when $A$ is understood. $L(W)$ is a $(j-1)$--dimensional tropical linear subspace and, being piecewise linear,
the subspace $L(W)$ has
a finite number of vertices  (however, an upper bound on how many  is not known in all cases; see \cite{Speyer}). Returning to the
question of which points are  vertices of $\p(A)$, the answer is that
the vertices 
of $\p(A)$ are all the vertices of all
subspaces $L(W)$, for  $W\in {{[d+1]}\choose{j}}$.  \label{para:vertices} The case $j=1$ gives
the $d+1$ generators of $\p(A)$.

\medskip
The easiest alcoved polytopes are \emph{boxes} and \emph{cubes}, determined by equations $x_i=cnst$. 
We fix a convenient matrix notation  for boxes with special matrices VNI and SNI (see Items \ref{item:a} and \ref{item:b} in
p.~\pageref{item:a}). Recall that  translation of an alcoved polyhedron $\p(A)$ is achieved by  conjugation of matrix $A$.
\begin{nota}[Box matrices]\label{nota:box_matrices} Given  real numbers $\ell_i>0$, $i\in [d]$, consider
\begin{enumerate}
\item  $B^{VNI}(d+1; \ell_1,\ell_2,\ldots,\ell_d)=[b_{i,j}]\in M_{d+1}(\R)$ with
$b_{i,j}=\begin{cases}
-\ell_i,&  d+1\neq i\neq j,\\
0,&\text{ otherwise,}\\
\end{cases}$. This matrix is VNI (easily checked)  and called the \emph{VNI box matrix with edge--lengths $\ell_j$.}
In particular, we have the \emph{VNI cube matrix}  $Q^{VNI}(d+1; \ell):=B^{VNI}(d+1; \ell,\ldots,\ell)$.

\item The conjugate matrix $D\odot B^{VNI}(d+1; \ell_1,\ell_2,\ldots,\ell_d)\odot D^{-1}$ is SNI (easily checked),
where $D=\diag(\ell_1/2,\ell_2/2,\ldots,\ell_d/2,0)$.
It is denoted $B^{SNI}(d+1; \ell_1,\ell_2,\ldots,\ell_d)=[c_{i,j}]$ and we have
$c_{i,j}=\begin{cases}
-\ell_i/2,&  j=d+1,\\
-\ell_j/2,&  i=d+1,\\
0,& i=j,\\
(-\ell_i-\ell_j)/2,& \text{otherwise.}\\
\end{cases}$
Similarly we have the cube matrix $Q^{SNI}(d+1; \ell)$.

\item A \emph{box matrix} is any diagonal conjugate of  the above,  i.e., $D'\odot B\odot {D'}^{-1}$,
where $D'=\diag(d'_1,d'_2,\ldots,d'_d,0)$ with  $d'_j\in \R$ and $B=B^{VNI}(d+1;\ell_1,\ell_2,\ldots,\ell_d)$.
It is NI (easily checked).
\end{enumerate}
\end{nota}

\begin{dfn}[from de la Puente \cite{Puente_QE}]\label{dfn:perturbation}
Any non--positive real matrix $E\in M_{d+1}(\R)$ with  null diagonal, last row and column is called \emph{perturbation matrix}.
In symbols, $E=[e_{i,j}]$ with $e_{i,i}=e_{d+1,i}=e_{i,d+1}=0$ and $e_{i,j}\le0, \forall i,j$.
\end{dfn}
In \cite{Puente_QE} it is proved that
for any  NI matrix $A\in M_{d+1}(\R)$  (not necessarily VNI or SNI), there exists a
unique decomposition $A=B-E$, where $B$ is a NI box matrix and   $E$ is a perturbation matrix.
The polytope  $\p(B)$ is called the \emph{bounding box} of the alcoved polytope $\p(A)$. \label{dfn:bounding_box} It is also proved
that $E$ is invariant under conjugation by diagonal matrices with zero last diagonal entry.

\section{Definition, characterization  and $f$--vector of IAPs}\label{sec:isocanted}
\begin{dfn}[Isocanted alcoved polytope (IAP)]\label{dfn:isocanted}
Let $A\in M_{d+1}(\R)$ be a NI matrix with decomposition $A=B-E$. The  alcoved polytope $\p(A)$  is \emph{isocanted} if $E$ is a \emph{constant perturbation matrix}, i.e., there exists $a>0$  such that $e_{i,j}=-a$, for all
$i,j\in[d]$, $i\neq j$. The number $a$ is called \emph{cant parameter} of $\p(A)$. We write $E=[-a]$, by abuse of notation.
\end{dfn}

\begin{rem}\label{rem:sym}
Every box in $\R^d$ is centrally symmetric and, by a translation, we can place its center of symmetry  at the origin of $\R^d$. An
IAP is a perturbed box with constant (whence symmetric) matrix $E$. Then,   every IAP  is centrally symmetric,by Item \ref{item:b} in p.~\pageref{item:b}.
\end{rem}

\begin{nota}[Special matrices for  visualized IAPs  and symmetric IAPs, with cubic bounding boxes]\label{nota:isocanted_matrices}
Given real numbers $a,\ell$,  consider the constant perturbation matrix  $E=[-a]\in M_{d+1}(\R)$ as above and the matrices (as in Notation \ref{nota:box_matrices})
\begin{enumerate}
\item $I^{VNI}(d+1; \ell,a):=Q^{VNI}(d+1;\ell)-E$,\label{item:11}
\item $I^{SNI}(d+1; \ell,a):=Q^{SNI}(d+1;\ell)-E$.
\end{enumerate}
It is an easy computation to check that, for these matrices  to be NI, it is necessary and sufficient that $0<a<\ell$.\footnote{The limit case $a=\ell$ provides a polytope of dimension less than $d$. The limit case $a=0$ provides the $d$--cube. Matrices $I^{VNI}(d+1; \ell_1,\ell_2,\ldots,\ell_d,a)$ and $I^{SNI}(d+1; \ell_1,\ell_2,\ldots,\ell_d,a)$ may be similarly defined, for $0<a<\min_j \ell_j$, but we will not use them.}
\end{nota}

The following is   the crucial step of the paper. Its proof contains the only tropical computations in what follows.

\begin{The}[Characterization of IAPs]\label{thm:charact}
An alcoved $d$--polytope $\PP=\p(A)$  is \emph{isocanted}  if and only if, for each    ${1}\le j\le d$ and  each
$W\in{ {[d+1]}\choose{j}}$ , the tropical linear subspace $L_A(W)$ has a unique vertex.
\end{The}
\proof
Without loss of generality, we can assume that the bounding box of $\PP$ is a cube (of  edge--length $\ell>0$) since an affine bijection  does not affect the result. We can also assume  that $\PP$ is located in $d$--space
so that ${\bf{0}}=\max\PP$, because a translation does not affect the result. Then $\PP=\p(C)$, with $C=Q^{VNI}(d+1;\ell)-E$,
for some positive  $\ell$, as in
in Notation \ref{nota:box_matrices} and Definition \ref{dfn:perturbation}. For $W\subset [d+1]$, let $C(W)$ denote the
$(d+1)\times j$ sized matrix  whose
columns are indexed by $W$ and taken from $C$.

$(\Rightarrow)$ Assume $\PP$ is IAP. Then  $E=[-a]$ is constant and then $C=I^{VNI}(d+1;\ell,a)=[c_{i,j}]$, as in
Item \ref{item:11} of Notation \ref{nota:isocanted_matrices}. In symbols,
$c_{i,j}=\begin{cases}
-\ell,&  i\neq j=d+1,\\
0,& i=j \text{ or } i=d+1,\\
-\ell+a,& \text{otherwise,}\\
\end{cases}$ \ with $0<a<\ell$.  Note that the tropical rank of $C$ is $d+1$ (meaning that the maximum in the
tropical permanent\footnote{The  \emph{tropical permanent} is the maximum of a collection of terms
(the definition mimics the classical one). Tropical permanent and \emph{tropical determinant} mean the same, in this paper. \emph{Tropical Laplace expansions} are one way to expand tropical determinants. For tropical permanent and
tropical rank issues, see  \cite{Butkovic_Libro,Develin_Santos_Sturm,Merlet}.} of $C$
is attained only once.\footnote{We have $\per_{tr}C=0$, attained only at the identity permutation.})
In particular,   $\rk_{\tr}C(W)=j$, for each proper subset $W\in {{[d+1]}\choose{j}}$.

For $j=1$,  $L(W)$ reduces to a point (a generator) and uniqueness is trivial.
Consider a point 
$x\in \R^{d+1}$ with $x_{d+1}=0$,   and let $C(W,x)$ be the matrix $C(W)$ extended
with column $x$. It is well--known (see  \cite{Richter,Speyer,Speyer_Sturm}) that $x\in L(W)$ if and only if $\rk_{\tr}C(W,x)\le j$,
(meaning that the maximum in \emph{each}  order $(j+1)$ tropical
minor\footnote{By \emph{tropical minor} we mean the tropical permanent (or determinant) of a square submatrix.
It is the maximum of a collection of terms.} is attained, at least, twice). Besides,
$x$ is a  vertex in $L(W)$ if and only if  the maximum in \emph{each} order $(j+1)$ tropical minor  of $C(W,x)$ is attained $(j+1)$
times. Indeed, the vertices of $L(W)$ are
got by computing
the corner locus of $L(W)$, then the corner locus of the corner locus, repeatedly.  Each iteration reduces the dimension of the
computed set, because  points where the maxima are attained one more time than previously, are computed.

\medskip
For each $2\le j\le d$ and each index family $1\le i_1<i_2<\cdots<i_j\le d+1$, let $m_{i_1,i_2,\ldots,i_j}$ (resp. $m_{i_1,i_2,\ldots,i_j}(x)$) denote the order $j$ minor of $C(W)$ (resp. $C(W,x)$) using rows $i_1,i_2,\ldots,i_j$.
Two cases arise.
    \begin{enumerate}
    \item \emph{If $d+1\notin W$,} then  it can be seen that  $m_{i_1,i_2,\ldots,i_j}=h(-\ell+a)$, where $h=|\{i_1,i_2,\ldots,i_j\}\setminus (W\cup\{d+1\})|$. In particular, $m_{i_1,i_2,\ldots,i_j}=0$, when $\{i_1,i_2,\ldots,i_j\}\subseteq W\cup\{d+1\}$.
    \item \emph{If $d+1\in W$,} then   $m_{i_1,i_2,\ldots,i_j}=h_1(-\ell)+h_2(-\ell+a)$, where
    $h_1=\begin{cases}
    1,& i_j\neq d+1,\\
    0,&\text{\ otherwise,}
    \end{cases}$ and $h_2=|\{i_1,i_2,\ldots,i_{j-1}\}\setminus W|$. In particular, $m_{i_1,i_2,\ldots,i_j}=0$, when $\{i_1,i_2,\ldots,i_j\}\subseteq W$.
    \end{enumerate}
 The  order $(j+1)$ minors in $C(W,x)$, expanded by the last column  by the tropical Laplace rule,  are
 \begin{equation}\label{eqn:minors_1}
 m_{i_1,i_2,\ldots,i_{j+1}}(x)=\max_{k\in[j+1]}\{x_{i_k}+m_{i_1,\ldots,i_{k-1},i_{k+1},\ldots, i_{j+1}}\}
 \end{equation}
 with $1\le i_1<i_2<\cdots<i_{j+1}\le d+1$,
 and the requirement that the maximum is attained $(j+1)$ times simply means that all the terms inside the maximum are equal, i.e.,
 \begin{equation}\label{eqn:minors_2}
 x_{i_k}+m_{i_1,\ldots,i_{k-1},i_{k+1},\ldots, i_{j+1}}=x_{i_{k'}}+m_{i_1,\ldots,i_{{k'}-1},i_{{k'}+1},\ldots, i_{j+1}},
 \ \forall k,k'\in[j+1].
 \end{equation}
  \begin{enumerate}
 \item \emph{If $d+1\notin W$,} then  $x_k=-\ell+a=c_{kj}$, for all $k\notin W\cup\{d+1\}$   (because $\rk_{\tr}C(W,x)\le j$ tells us that $x$ is a tropical affine combination of the columns in $C(W)$),  and equalities  (\ref{eqn:minors_2}) imply $x_k=0$, for all $k\in W\cup\{d+1\}$.
     \item \emph{If $d+1\in W$,} then  $x_k=x_{k'}$, for all $k,k'\notin W$ (because $\rk_{\tr}C(W,x)\le j$ tells us that $x$ is a tropical affine combination of the columns in $C(W)$), and equalities  (\ref{eqn:minors_2}) imply $x_k=-a$, for all $d+1\neq k\in W$, $x_k=-l$, for all $k\notin W$, $x_{d+1}=0$.
\end{enumerate}

\medskip

$(\Leftarrow)$ We have $\PP=\p(C)$, where $C=Q^{VNI}(d+1;\ell)-E$  is a NI matrix. Assume that, for each ${1}\le j\le d$ and each
$W\in {{[d+1]}\choose j}$, the
tropical linear subspace $L(W)$ has a unique vertex denoted $x_W^*$. We write $x^*$ whenever $W$ is understood. We have $x^*_{d+1}=0$.

 Since  $\PP$ is centrally symmetric, then by Remark \ref{rem:sym}, the matrix  $E$ is symmetric. We want to prove that
$E$ is constant.
Fix $w\in[d]$ and take $W=\{w, d+1\}$. Use that for each  order 3   minor of the matrix $C(W,x^*)$ (where $x^*$ depends on $w$) all terms in the maximum
are equal.
Considering those minors involving  three different indices $i,w,d+1$, we get
\begin{equation}
x^*_i+m_{w, d+1}=x^*_w+m_{i, d+1}=x^*_{d+1}+\begin{cases}
m_{i,w}&\text{if\ } i<w\\
m_{w,i}&\text{otherwise\ }
\end{cases}=0-\ell
\end{equation}
whence
\begin{equation}\label{eqn:iw}
x^*_i+0=x^*_w-\ell-e_{i,w}=-\ell,
\end{equation}
and so $e_{w,i}=e_{i,w}=x^*_w$. Letting $i\in [d]$ vary in (\ref{eqn:iw}), we get that $E=[x^*_w]$ is constant.
\endproof

\begin{ex}\label{ex:proof_thm}
Let $d=5$.  If $j=5$ and  $W=[5]$, then the tropical Laplace expansion by the last column yields  $\per_{tr}C(W,x)=\max\{x_1,x_2,x_3,x_4,x_5,0\}$. This maximum is attained  by all terms if and only if $x_k=0$, all $k\in [5]$, so the unique vertex of $L(W)$
is the origin.

 If  $j=3$ and  $W=[3]$ then
$C(W,x)=\left[\begin{array}{rrrr}0&-\ell+a&-\ell+a&x_1\\-\ell+a&0&-\ell+a&x_2\\-\ell+a&-\ell+a&0&x_3\\-\ell+a&-\ell+a&-\ell+a&x_4\\
-\ell+a&-\ell+a&-\ell+a&x_5\\0&0&0&0\end{array}\right]$.  Since $x$ is a tropical affine combination of the columns of $C(W)$, it follows that $x_4=x_5=-\ell+a$. Since the maximum
\begin{equation}
m_{1234}(x)=\max\{x_1+m_{234},x_2+m_{134},x_3+m_{124},x_4+m_{123}\}=\max\{x_1-\ell+a,x_2-\ell+a,x_3-\ell+a,x_4\}
\end{equation}
is attained by all terms, we get
\begin{equation}
x_1-\ell+a=x_2-\ell+a=x_3-\ell+a=x_4=-\ell+a
\end{equation}
whence $x_1=x_2=x_3=0$. The unique vertex of $L(W)$ is $[0,0,0,-\ell+a,-\ell+a]^T$.

If $j=3$ and $W=\{1,2,d+1\}$ then $C(W,x)=\left[\begin{array}{rrrr}0&-\ell+a&-\ell&x_1\\-\ell+a&0&-\ell&x_2\\-\ell+a&-\ell+a&-\ell&x_3\\-\ell+a&-\ell+a&-\ell&x_4\
\\-\ell+a&-\ell+a&-\ell&x_5\\0&0&0&0\end{array}\right]$.  Since $x$ is a tropical affine combination of the columns of $C(W)$, it follows that $x_3=x_4=x_5$. Since the maximum
\begin{equation}
m_{1236}(x)=\max\{x_1+m_{236},x_2+m_{136},x_3+m_{126},m_{123}\}=\max\{x_1-\ell+a,x_2-\ell+a,x_3,-\ell\}
\end{equation}
is attained by all terms, we get
\begin{equation}
x_1-\ell+a=x_2-\ell+a=x_3=-\ell
\end{equation}
whence $x_1=x_2=-a$ and $x_3=x_4=x_5=-\ell$. The unique vertex of $L(W)$ is $[-a,-a-\ell,-\ell,-\ell]^T$.
\end{ex}

\begin{rem}
We have   $x_W^*=\bigoplus_{j\in W}\underline{j}$, whenever $d+1\notin W$.
\end{rem}

\begin{cor}
[Bijection on set of  vertices of IAP]\label{cor:vertices}
Given any isocanted alcoved $d$--polytope $\PP$,
the vertices of $\PP$ are in bijection with the \emph{proper}  subsets $W\subset [d+1]$.

\end{cor}
\proof
As a a set, a tropical line is  a finite union of  classical segments and
halflines.\footnote{A balance condition at each point of each tropical algebraic variety
is satisfied, but we do not use it in this paper.}
As a set, a tropical segment\footnote{A tropical segment is the family of all tropical affine combinations of  two points.} is a
finite union of  classical segments.
The tropical line strictly contains the tropical segment determined by two given points, and the difference set
is a finite union of halflines; see \cite{Develin_Sturm,Puente_line,Richter}.
For an alcoved polytope $\PP$, this implies that the \emph{skeleton}\footnote{The skeleton is the 1--dimensional subcomplex of the border
complex $\partial \PP$. It is a graph, whose diameter is computed in Corollary \ref{cor:diameter}.} of $\PP$ is \label{dfn:skeleton}
contained in  the 1--dimensional complex $\bigcup_{W\in {[d+1]\choose2}}L(W)$.  For each $W$
with $|W|=2$, the set $L(W)\setminus \PP$ is a finite union of halflines.
Every \emph{generated} vertex of $\PP$ is also a vertex of the complex  $\bigcup_{W\in {{[d+1]}\choose2}}L(W)$, and every edge of $\PP$
is contained in an edge of  $\bigcup_{W\in {{[d+1]}\choose2}}L(W)$. The containment is strict  exactly for those edges of $\PP$ emanating
from generators.

If $\PP$ is IAP and $\underline{i},\underline{j}$ are two generators (with $i,j\in[d+1], i<j$),
the tropical line determined by them has a unique vertex, which will be denoted $\underline{ij}$. If $i,j,k\in[d+1]$ with $i<j<k$,
the tropical plane determined by them has a unique vertex, which will be denoted $\underline{ijk}$. It can be checked
that $\underline{ijk}$ is the unique vertex of the tropical line determined by $\underline{ij}$ and $\underline{k}$. Recursively,
vertices of $\PP$ are labeled in this fashion. The stated bijection follows.
\endproof

\begin{nota} The label of
the vertex corresponding to $W\subset [d+1]$ is $\underline{W}$ (underlined). The  cardinality $|W|$ is called
\emph{length} of $\underline{W}$.
\end{nota}

\begin{nota}[Parent and child]\label{nota:edges}
Assume $\PP$ is an isocanted  alcoved $d$--polytope. Two vertices  in $\PP$ are joined by an edge in $\PP$ if and only if they are  labeled
$\underline{W}$ and $ \underline{W'}\subset [d+1]$ with  $\emptyset\neq W\subset W'$ and $|W|+1=|W'|$.  We say that
$\underline{W}$ is a \emph{parent}
of $\underline{W}'$ and  $\underline{W'}$ is a \emph{child} of $\underline{W}$. A 2--face of $\PP$ is determined by four
vertices with labels $\underline{jW}$, $\underline{jkW}$, $\underline{jrW}$, $\underline{jkrW}$, with
$W\subset[d+1]\setminus\{j,k,r\}$, for $j,k,r\in[d+1]$ pairwise different.\footnote{$jW$ is shorthand for $\{j\}\cup W$.}
\end{nota}

\begin{The}[$f$--vector for IAP]\label{thm:f_vector_IAP}
$I_{d,j}=(2^{d+1-j}-2){{d+1}\choose{j}}, \quad 0\le j\le d-1$.
\end{The}
\proof

First, $I_{d,0}=\left|\bigcup_{j=1}^d{{[d+1]}\choose{j}}\right|=2^{d+1}-2$ is the number of proper subsets of $[d+1]$.

Second, the number of facets is   $I_{d,d-1}=(d+1)d$ by (\ref{eqn:ineq_ij}) and (\ref{eqn:ineq_i}). Another proof is this:  as
we mentioned in
p.~\pageref{beveling}, an alcoved polytope is obtained from a box, where  we may cant only
the $(d-2)$--faces not meeting  two distinguished opposite vertices; thus, we may cant  half of the $(d-2)$--faces of the box. In
an IAP we do cant every cantable $(d-2)$--face, therefore $I_{d,d-1}=B_{d,d-1}+B_{d,d-2}/2=(d+1)d$.

\medskip
For $1\le j\le d$, the number of vertices of length $j$  is ${d+1}\choose{j}$, by Theorem \ref{thm:charact}.

\medskip
Assume $2\le j\le d$. A vertex of length $j$ has  $j$ parents, by Notation \ref{nota:edges}. The total number of edges
is $\sum_{j=2}^{d}{{d+1}\choose{j}}j=(d+1)\sum_{j=2}^{d}{d\choose{j-1}}=(d+1)\sum_{k=1}^{d-1}{d\choose{k}}=(d+1)(2^{d}-2)=I_{d,1}$, (where we have used the
equalities  ${{d+1}\choose{j}}j=(d+1){{d}\choose{j-1}}$ and $2^d=\sum_{j=0}^{d} {d \choose j}$).

Assume $3\le j\le d$.  A vertex of length $j$ has  ${{j}\choose{2}}$ grandparents (i.e., parent of parent). The total number of 2--faces
is $\sum_{j=3}^{d}{{d+1}\choose{j}}{{j}\choose{2}}={{d+1}\choose{2}}\sum_{j=3}^{d}{{d-1}\choose{j-2}}
={{d+1}\choose{2}}\sum_{k=1}^{d-2}{{d-1}\choose{k}}={{d+1}\choose{2}}(2^{d-1}-2)=I_{d,2}$
(where we have used the equality ${{d+1}\choose{j}}{{j}\choose{2}}={{d+1}\choose{2}}{{d-1}\choose{j-2}}$).

Similarly, the total number of $r$--faces
is $\sum_{j=r+1}^{d}{{d+1}\choose{j}}{{j}\choose{r}}={{d+1}\choose{r}}\sum_{j=r+1}^{d}{{d+1-r}\choose{j-r}}
={{d+1}\choose{r}}\sum_{k=1}^{d-r}{{d+1-r}\choose{k}}={{d+1}\choose{r}}(2^{d+1-r}-2)=I_{d,r}$
(where we have used the equality ${{d+1}\choose{j}}{{j}\choose{r}}={{d+1}\choose{r}}{{d+1-r}\choose{j-r}}$).
\endproof

\begin{rem}
A $d$--IAP is a canted box where \emph{all} cantable $(d-2)$--faces are canted. On the contrary, alcoved polytopes exist where some cantable $(d-2)$--faces of the bounding box  remain uncanted. Among alcoved polytopes, IAPs are  maximal in facets because in an IAP  we cant every possible cantable $(d-2)$--face. Notice that IAPs are neither simplicial nor simple and far from being neighborly.
\end{rem}

\begin{rem}
Notice the coincidence of $I_{d,j}$ with  the triangular sequence OEIS A259569
(collecting  the number of $j$--dimensional faces on the polytope that is the convex hull of all permutations of the
list $(0,1,\ldots,1,2)$, where there are $d - 1$ ones). Also notice the coincidence of $I_{d,j}$ with  the absolute values of
the triangular sequence  OEIS A138106 (collecting the coefficients of the Taylor expansion around the origin of
the function of two variables $p(x,t)=e^{(x-2)t} - 2e^{(x-1)t}$. Functions of  similar appearance are called Morse potentials);
see \cite{OEIS}.
\end{rem}

\medskip
The study of \emph{cubical polytopes} began in the late 1990's in \cite{Blind_fewer, Blind_almost_simple}. \emph{Zonohedra} were  first considered  by the crystallographer E.S. Fedorov, by the end of the XIX century.
In the rest of this section, we  prove that IAPs are cubical polytopes and zonohedra.

A $d$--\emph{cuboid} is a polytope combinatorially equivalent to a $d$--cube.
A $d$--cuboid is denoted $\CC^d$. A polytope is \emph{cubical}\label{dfn:cubical}
if every face in it is a cuboid (equivalently,
if every facet in it is a cuboid).
A $d$--polytope is \emph{almost simple} if the valence of each vertex is $d$ or $d+1$. A
$d$--polytope  $\PP$
is \emph{liftable (to a $(d+1)$--cuboid)} if its boundary complex $\partial  \PP$ is combinatorially equivalent to a subcomplex of
the complex
$\partial  \CC^{d+1}$.

Take any vertex $V$ in a cuboid
$\CC^{d+1}$ and consider the subcomplex $\FF^d_V$ of $\partial  \CC^{d+1}$ determined by the facets of $\CC^{d+1}$ meeting $V$.
Consider the
subcomplex $\CC^d_V$ of $\FF^d_V$  determined by the outer faces of
$\FF^d_V$ (the underlying set of  $\CC^d_V$ is  $\partial  \FF^{d}_V$).
A 
polytope $\PP$ is \emph{$d$--elementary} if the complex $\partial  \PP$
is combinatorially equivalent to the subcomplex  $\CC^d_V$. We call $\FF^{d}_V$ \emph{(cuboid) cask at $V$}.\label{dfn:cube_cask}
\label{dfn:elementary}

It is clear that $d$--elementary is more specific than liftable.
Saying that  $\PP$ is  $d$--elementary
means that $\PP$ is   (combinatorially equivalent to)  the pasting of $d+1$ $d$--cuboids all having a vertex $V$ in common. More generally,   $k$--elementariness describes the property of $\PP$ being combinatorially equivalent to the pasting
of $k+1$ $d$--cuboids, all sharing a $(d-k)$--face.   In particular,
a $d$--cuboid is 0--elementary.  A $k$--elementary $d$--polytope   is
obtained from a $(k-1)$--elementary polytope by pasting (combinatorially) a $d$--cuboid to it. A $k$--elementary $d$--polytope is denoted $\CC^d_k$.

The main theorem in \cite{Blind_fewer} states that if $d\ge4$
and $\PP$ is a cubical $d$--polytope, then $\PP$ is $k$--elementary, for some $k$ with $0\le k\le d$. It is also proved that both
$\CC^d_{d-1}$ and $\CC^d_{d}$ have $2^{d+1}-2$ vertices, while $\CC^d_k$ has fewer than $2^{d+1}-2$ vertices, for other values of $k$.
Corollary 1 in \cite{Blind_almost_simple} states that,  for $d\ge4$, a $d$--polytope is liftable if and only if  it is cubical,
almost simple and has, at most, $2^{d+1}$ vertices.

\begin{cor}\label{cor:face_lattice;unique;cubical} For each $d\ge2$,
\begin{enumerate}
\item the face lattice of a $d$--IAP is the lattice of proper subsets of $[d+1]$,\label{item:lattice}
\item on  the set of vertices of a $d$--IAP, the mapping $\underline{W}$ to $\underline{[d+1]\setminus W}$ is an involution, \label{item:duality}
\item there exists a unique combinatorial type of $d$--IAP,
\item every IAP is cubical and almost simple. \qed
\end{enumerate}
\end{cor}
\proof
\begin{enumerate}
\item This is direct consequence of Corollary \ref{cor:vertices}.
\item This is  due to the lattice order--reversing isomorphism $W\mapsto [d+1]\setminus W$.
\item This is immediate from  Item \ref{item:lattice}.
\item Let $\PP$ be an IAP and  $\BB$ be the bounding box of $\PP$ (defined in p.~\pageref{dfn:bounding_box}). The $(d-2)$--faces meeting the  two distinguished vertices of $\BB$ are not cantable and so, those two points are vertices of both  $\PP$ and $\BB$, and they have  the same valence in  $\PP$ and $\BB$ (the valence   is $d$ in $\BB$).  In $\PP$ one of these vertices is the generator $\underline{d+1}$ and the other one has label
    $\underline{12\ldots d}$. All generators (resp. $d$--generated vertices) of $\PP$ have the same valence.  Generators do not have parents and vertices of length $d$ do not have children. Now, for $2\le j\le d-1$, the  valence of a vertex of $\PP$ of length $j$ is the sum of the number of parents and number of children, namely, $j+(d+1-j)=d+1$.
\end{enumerate}
\endproof
For $d=2$, an IAP is a hexagon (with slopes 0,1,$\infty$) and every vertex in it has valence 2.  For $d=3$, an IAP is  combinatorially equivalent to a
\emph{rhombic dodecahedron}, whose $f$--vector is $(14,24,12)$.

\begin{nota}
Since the combinatorial type is unique,  we can fix a notation for a  $d$--IAP: it  is denoted $\II^d$ in what follows.
\end{nota}
\begin{cor}\label{cor:element_lift}
$\II^d$ is  $d$--elementary, for $d\ge2$.
\end{cor}
\proof
For $d\ge4$, $k$--elementariness  follows from the  main theorem in \cite{Blind_fewer}, and
$f_0(\CC^d_{d-1})=f_0(\CC^d_d)=f_0(\II^d)$ tell us that $k=d-1$ or $d$.  The generator $\underline{d+1}$ (also the vertex $\underline{12\ldots d}$) plays the role of vertex $V$ in the definition above in p.~\pageref{dfn:elementary}, so that $k=d$ follows.

$2$--elementariness is   easy for $d=2$:  $\II^2$ is a  hexagon, and it is combinatorially equivalent to $\CC_V^{2}$,
which is the border complex of
a cube \emph{cask} $\FF^{2}_V$ at a vertex $V$ of the cube. For $d=3$, extended explanations are given in section \ref{sec:d=3,4}.
\endproof

The $f$--vector  of a cask $\FF^d_V$ clearly is
\begin{equation}\label{eqn:f_vector_box_cask}
C_{d,j}=\left(2^{d-j}-1\right){{d}\choose{j}},\quad j=0,1,\ldots,d-2.
\end{equation}
Since  $\II^d$ is $d$--elementary, then (\ref{eqn:f_vector_box_cask}) and (\ref{eqn:Idj}) satisfy the following the relation
\begin{equation}\label{eqn:relation}
I_{d,j}=2C_{d,j}+I_{d-1,j-1},
\end{equation}
which has the practical application that, in order to understand $\partial\II^d$, it is enough that we  look a two  cube casks and
one belt joining them. See
section \ref{sec:d=3,4} for details in dimensions 3 and 4.

\medskip
Recall that a \emph{zonotope}  is a (Minkowski) sum of segments. A known characterization of zonotope is that it is a polytope all whose 2--faces are centrally symmetric (see \cite{chap15}),
and this is satisfied by IAPs. A direct proof is given below.

\begin{cor}\label{cor:zono}
Every IAP is a zonotope.
\end{cor}
\proof $\II^d$ is obtained from a $d$--box
$\BB=\BB(\ell_1,\ell_2,\ldots,\ell_d)\subset \R^d$ with $\max \BB$ at the origin,
edge--lengths $\ell_j>0$ and cant parameter $a$ with $0<a<\min_{j} \ell_j$, and  $\II^d=\BB+[0,a v_{d+1}]$ holds true,
where $(v_1,v_2,\ldots,v_d)$ is the standard basis in $\R^d$ and $v_{d+1}=v_1+v_2+\cdots+v_d$.
\endproof

\section{Cases  $d=3$ and $4$.}\label{sec:d=3,4}
In this section we describe IAPs in the small dimensions, for   a better understanding of results proved in the previous section.
In addition, for $d=4$, we compute two well--known invariants (fatness and $f_{03}$).

Fix $d\ge2$.  Two opposite vertices in $\II^d$  are distinguished: $\NN:=\max \II^d$, called the
\emph{North Pole},
and $\SS:=\min \II^d$ called the \emph{South Pole}  of $\II^d$.\footnote{This idea, which  goes back to Kepler, has been developed for alcoved polytopes in  \cite{Puente_QE}.} The label
of $\NN$ is $\underline{1 2\cdots d}$, and the label of $\SS$ is $\underline{d+1}$ ($\SS$ is a generator).
The cask $\FF^{d}_{\NN}\subset \partial \II^d$  introduced in p.~\pageref{dfn:cube_cask} (resp. $\FF^{d}_{\SS}$) is called \emph{North Polar Cask}
(resp. \emph{South Polar Cask}) of $\II^d$.
Vertices included in the North (resp. South) Polar Cask are exactly those omitting (resp. including) digit
$d+1$ in their label.  The \emph{Equatorial Belt} is, by definition,  the subcomplex of $\partial \II^d$ determined by all faces of $\II^d$
not meeting the poles. The Equatorial Belt is the complex of all facets of $\II^d$ containing edges  in the  direction of vector $v_{d+1}=(1,1,\ldots,1)^T$. These are the edges  joining
vertices $\underline{W}$ and $\underline{W d+1}$,
for proper subsets  $W\subset [d]$. The complex $\partial \II^d$ is the union of the Polar Casks and the Equatorial Belt.

\medskip
A Polar Cask is homeomorphic to a closed $(d-1)$--disk.
The Equatorial Belt is homeomorphic to a closed $(d-1)$--cylinder, i.e.,  $S^{d-2}\times[-1,1]$ (the Cartesian product of a
$(d-2)$--sphere and a closed interval).

\begin{figure}[ht]
\centering
\includegraphics[width=10cm]{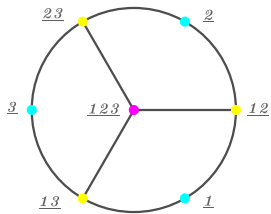}
\caption{North Polar Cask for $d=3$.}\label{fig1_NorthPolarCask_d=3}
\end{figure}

\begin{figure}[ht]
\centering
\includegraphics[width=10cm]{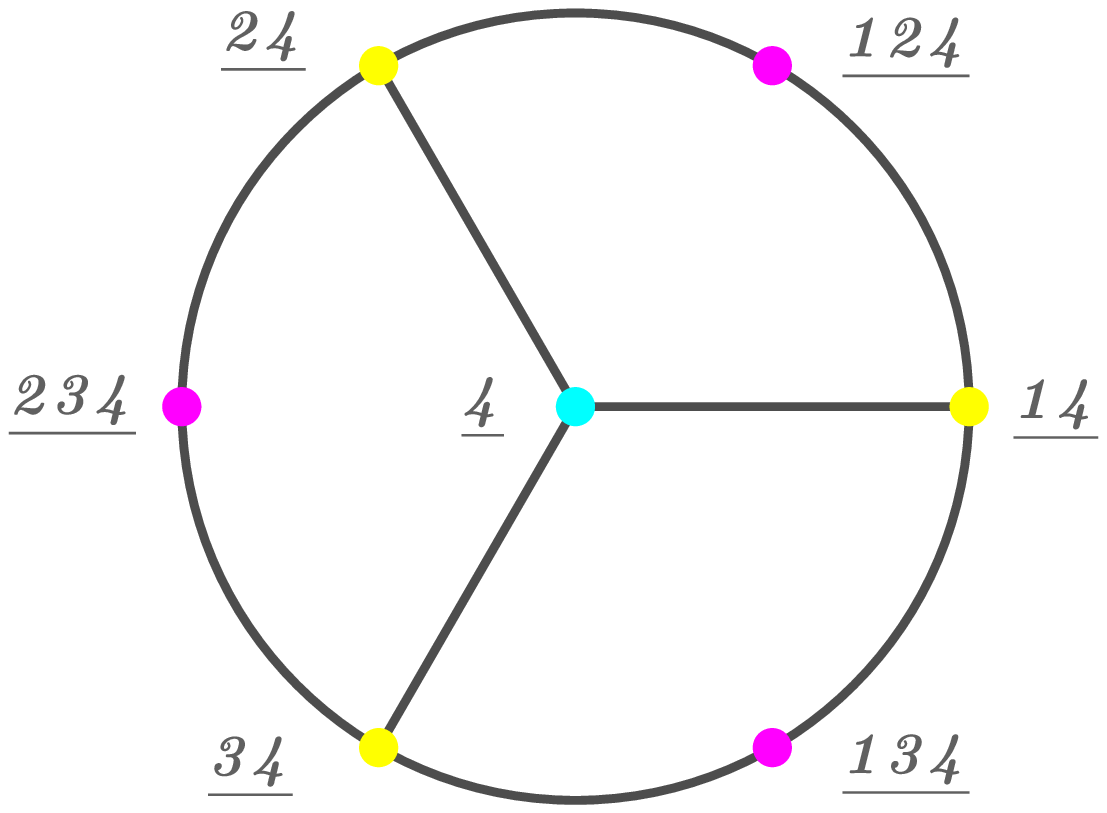}
\caption{South Polar Cask for $d=3$.}\label{fig2_SouthPolarCask_d=3}
\end{figure}

\begin{figure}[ht]
\centering
\includegraphics[width=6cm]{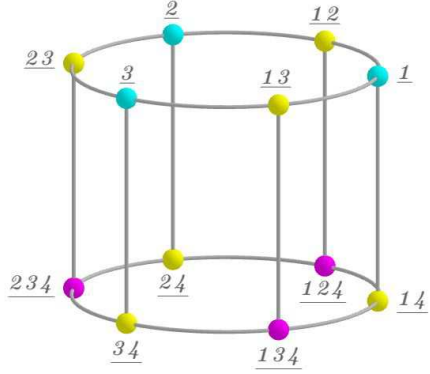}
\caption{Ecuatorial  Belt for $d=3$.}\label{fig3_EcuatorialBelt_d=3}
\end{figure}

\medskip
\textbf{Case $d=3$:} we have $\NN=\underline{123}$ and the North Cask is homeomorphic to a 2--disk with one interior point labeled $\underline{123}$,
points in the circumference labeled $\underline{1},\underline{12},\underline{2},\underline{23},\underline{3},\underline{13}$ and inner
edges joining $\underline{12},\underline{23},\underline{13}$ to $\underline{123}$
(see figure \ref{fig1_NorthPolarCask_d=3}). The South Pole is
$\SS=\underline{4}$ and the South Cask is homeomorphic to a 2--disk with one interior point labeled $\underline{4}$, points in the circumference
labeled $\underline{14},\underline{124},\underline{24},\underline{234},\underline{34},\underline{134}$ and inner edges
joining $\underline{14},\underline{24},\underline{34}$ to $\underline{4}$  (see figure \ref{fig2_SouthPolarCask_d=3}). The Equatorial Belt is homeomorphic to a cylindrical  surface (see figure \ref{fig3_EcuatorialBelt_d=3}).
Identification of borders of polar casks with border components of the
cylinder is easily done by using vertex labels. The $f$--vector of a $2$--polar cask is the sum of the $f$--vector of the circumference complex $(6,6)$ and of the internal subdivision $(1,3)$, yielding $(7,9)$, which agrees with $(C_{3,0},C_{3,1})$ in (\ref{eqn:f_vector_box_cask}).

\begin{figure}[ht]
\centering
\includegraphics[width=6cm]{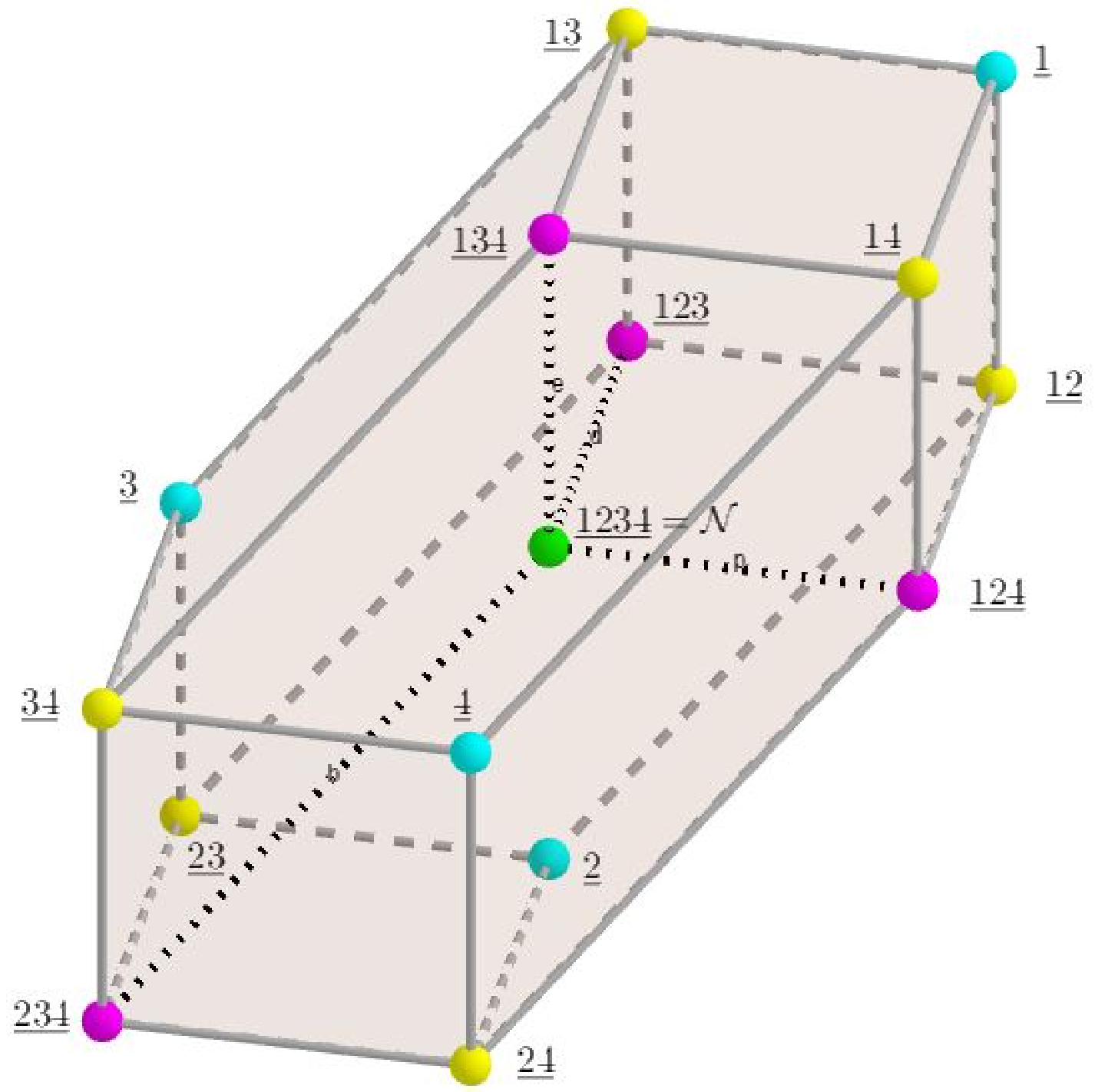}
\caption{North Polar Cask for $d=4$.}\label{fig4_NorthPolarCask_d=4}
\end{figure}

\textbf{Case $d=4$:} the North Cask is homeomorphic to a solid 3--sphere with one interior point
labeled $\NN=\underline{1234}$,
points on the surface labeled $\underline{i}$, $\underline{ij}$,  and $\underline{ijk}$, with $i,j,k\in[4]$, pairwise different.
Edges  join parent and child (see  Notation \ref{nota:edges}). Combinatorially, the cask is equivalent to a solid rhombic dodecahedron with an
interior point labeled $\underline{1234}$ and six quadrangular inner 2--faces given by  $\underline{ij}$, $\underline{ijk}$,
$\underline{ijr}$, $\underline{1234}$, with $\{i,j,k,r\}=[4]$ (see figure \ref{fig4_NorthPolarCask_d=4}).

\begin{figure}[ht]
\centering
\includegraphics[width=6cm]{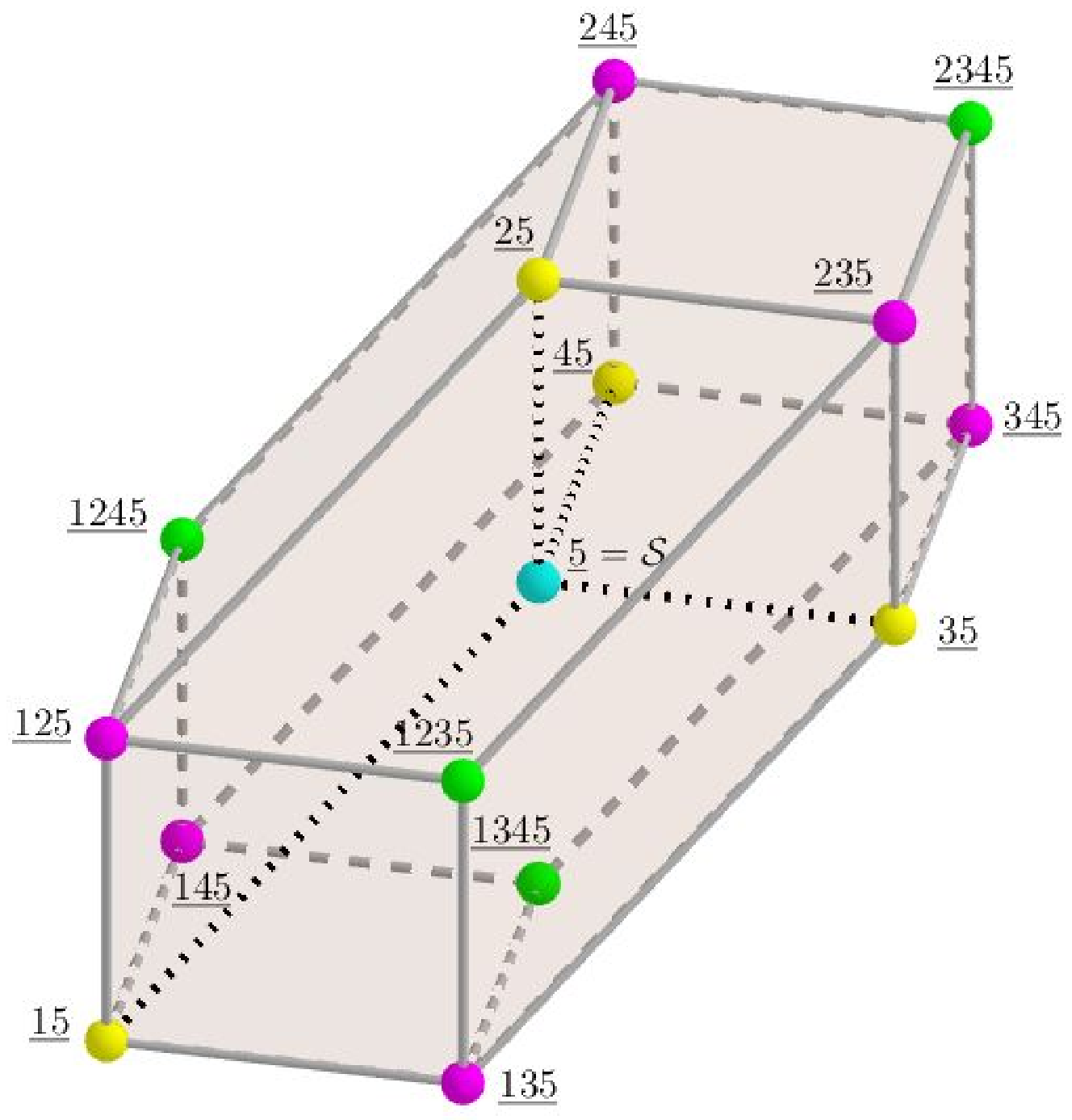}
\caption{South Polar Cask for $d=4$.}\label{fig5_SouthPolarCask_d=4}
\end{figure}

The South Cask is homeomorphic to a solid 3--sphere with one interior point
labeled $\SS=\underline{5}$, points on the surface labeled $\underline{i5}$, $\underline{ij5}$,  and $\underline{ijk5}$, with
$i,j,k\in[4]$, pairwise different.
Edges are determined  by Notation \ref{nota:edges}. Combinatorially, the cask is equivalent to a solid rhombic dodecahedron with an
interior point labeled $\underline{5}$ and six quadrangular inner 2--faces given by  $\underline{i5}$, $\underline{ij5}$,
$\underline{ik5}$, $\underline{ijk5}$, with $i,j,k\in [4]$ pairwise different (see figure \ref{fig5_SouthPolarCask_d=4}).
The $f$--vector of a rhombic dodecahedron
is $(14,24,12)$ and the internal subdivision adds $(1,4,6)$, so that the sum $(15,28,18)$ is the $f$--vector of  a
$3$--polar cask, which agrees with $(C_{4,0},C_{4,1},C_{4,2})$ in (\ref{eqn:f_vector_box_cask}).
The Equatorial Belt is homeomorphic to a 3--cylinder $S^2\times[-1,1]$. Identification of borders of polar casks with
border components of cylinder is easily done by using vertex labels.

\bigskip
Researchers are deeply interested in 4--polytopes, due to the peculiar properties they show (from the classification of the regular
ones obtained by Schl\"{a}fli in the XIX century,  to the  Richter-Gebert's Universality Theorem of 1996, which roughly says that the realization space of a 4--polytope can
be \lq\lq arbitrarily wild or ugly\lq\lq, see \cite{chap15}). Fatness is a convenient function to study the family ${\bf{F}}_4\subset \N^4$ of  $f$--vectors
of 4--polytopes. The  set ${\bf{F}}_4$ is   not well understood. The \emph{fatness} of a 4--polytope $\PP$ is defined as $F(\PP)=\frac{f_1+f_2-20}{f_0+f_3-10}$.
It is known that $F(\PP)\in [\frac{5}{2},3)$, for all simplicial and all simple $\PP$. It is also known that
$F(\PP)\le 5$, for all 4--zonotopes $\PP$ (see  \cite{Ziegler}). According to Ziegler,
\lq\lq the existence/construction of 4--polytopes of high fatness\rq\rq (greater than or equal to 9) \lq\lq is a
key problem.\rq\rq

$f$--vectors have been generalized in a number of ways. Generalizations considered in this paper are:
to count \emph{vertex--facet incidences} (denoted  $f_{03}$ below), to count \emph{flags}
(see Corollary \ref{cor:flag_conj}) and  the  \emph{cubical $g$--vector} (see Proposition \ref{prop:CLBC}).

\begin{rem} We have $I_4=(30,70,60,20)$ and \begin{enumerate}
 \item \emph{fatness} of $\II^4$ is $\frac{f_1+f_2-20}{f_0+f_3-10}=\frac{11}{4}$,
\item in $\II^4$ we have $f_{03}=160$ (since there are $I_{4,3}=(d+1)d=20$ 3--cubes
(with 8 vertices each) and no other 3--faces).
\end{enumerate}
\end{rem}
So fatness of IAPs will not surprise Ziegler!

\textbf{Key to colors}: blue dots are generators, yellow dots are vertices of length 2,  magenta dots are vertices of length 3, green dots are vertices of length 4.


\section{Five conjectures proved for IAPs}\label{sec:conjectures}
Consider the set $\MM$ of lower triangular infinite matrices with both entries \textbf{and} indices  in $\Z_{\ge0}$.
Examples of matrices in $\MM$ are the \emph{2--power matrix}, denoted $T$, defined by $T_{d,k}=\begin{cases}
2^{d-k},& 0\le k\le d,\\
0,&\text{otherwise},
\end{cases}$
and the \emph{Pascal matrix}, denoted $P$, defined by $P_{d,k}=\begin{cases}
{d\choose k},& 0\le k\le d,\\
0,&\text{otherwise.}
\end{cases}$ With the Hadamard or entry--wise product, multiply the former matrices, obtaining  $B:=T\circ P=P\circ T\in \MM$ and notice that the $d$--th row of
$B$ shows the $f$--vector of a $d$--box (padded with zeros),  for $d\in \Z_{\ge0}$; see (\ref{eqn:f_vector_box}). We call $B$ is the $f$--vector \emph{box matrix}.
Next, consider the matrix $H\in\MM$ defined by
\begin{equation}\label{eqn:Hdk}
H_{d,k}=\begin{cases}
(2^{d-k}-1){{d+1}\choose k},&0\le k\le d-1,\\
1/2,&k=d,\\
0,&\text{otherwise.}
\end{cases}
\end{equation}

For  fixed $d\ge2$, we  study the growth\footnote{$H_{d,k}$ is  an expression involving
2--powers and binomial coefficients.
Precisely,  $H_{d,k}=(T_{d,k}-1)P_{d+1,k}$ is the product of two factors. For
sufficiently  small  $k$,  the first factor dominates (meaning, is larger than the other factor), as  in the cases
$H_{d,0}=2^{d}-1$, $H_{d,1}=(2^{d-1}-1)(d+1)$ and  $H_{d,2}=(2^{d-2}-1)(d+1)d/2$. However, for
sufficiently  large  $k$, the second factor dominates, as in the cases
$H_{d,d-3}=7(d+1)d(d-1)(d-2)/24$, $H_{d,d-2}=(d+1)d(d-1)/2$ and $H_{d,d-1}=(d+1)d/2$.} of the sequence $H_{d,k}$,  with $0\le k< k+1\le d-1$.

\begin{prop}\label{prop:extremes}
For each $d\ge0$, we have $H_{d,d-1}\le H_{d,0}$ with equality only for $d=0,1,2$.
\end{prop}
\proof
The inequality $(d+1)d/2\le 2^d-1$ is easily proved by induction on $d$
(degree 2 polynomials grow slower than 2--powers.)
\endproof

Recall that a sequence $a_k$ is \emph{log--concave} if $a_{k+1}^2\ge a_{k}a_{k+2}, \forall k$; see \cite{Brenti,Stanley}.

\begin{prop}\label{prop:log_concave}
For $d\ge2$, the sequence  $\{H_{d,k}: 0\le k\le d-1\}$ is log--concave.
\end{prop}
\proof
For fixed  $d$, the sequence
$T_{d,k}-1=2^{d-k}-1$ is log--concave, because $(T_{d,k+1}-1)^{2}-(T_{d,k}-1)(T_{d,k+2}-1)=2^{d-k-2}>2>0$, for $0\le k\le d-3$.
It is easy to check that any row of
Pascal's triangle is a log--concave sequence. Since
 the  termwise product of two log--concave sequences (with the same number of
terms) is log--concave, then the result follows  for $H_{d,k}$.
\endproof
Notice $I_{d,k}=2 H_{d,k}$, for $0\le k{\le}d$.

\begin{cor}[Unimodality holds for isocanted]
For each $d\ge2$, the sequence  $\{I_{d,k}: 0\le k\le d-1\}$ is unimodal.
\end{cor}
\proof
It is easy to show that  every log--concave sequence is unimodal (but not conversely). The sequence $H_{d,k}$ is unimodal and so is its double.
\endproof

\begin{prop}\label{prop:d_tercios}
For fixed $d\ge2$, the maximum in the sequence $I_{d,k}$ is attained at the integer $\lfloor\frac{d}{3}\rfloor$.
\end{prop}
\proof Cases $d=2$, 3 and 4 are checked directly (the $f$--vectors are $(6,6)$, $(14,24,12)$ and $(30,70,60,20)$). Assume
$d\ge5$  and $0\le k\le d-2$.  Define the quotient
\begin{equation}\label{eqn:quotient}
Q_{d,k+1}:=\frac{I_{d,k+1}}{I_{d,k}}=\frac{(2^{d-k-1}-1)(d-k+1)}{(2^{d-k}-1)(k+1)}
\end{equation}
and the terms
\begin{equation}\label{eqn:equality}
L_{d,k+1}:=2^{d-k-1}(d-3k-1), \quad R_{d,k+1}:=d-2k.
\end{equation}
We have $I_{d,k+1}\ge I_{d,k}$ if and only if  $Q_{d,k+1}\ge1$ if and only if $L_{d,k+1}\ge R_{d,k+1}$, because we have
cleared the positive denominator in (\ref{eqn:quotient}) and grouped terms.
 The exponent $d-k-1$  appearing in $L_{d,k+1}$ is at least 1.
The sign of the factor $(d-3k-1)$ in $L_{d,k+1}$ is not constant. We have
$\frac{d}{3}\le\frac{2d-5}{3}$, since $d\ge5$.
We prove
\begin{enumerate}
\item  if $k\le \frac{d-2}{4}$, then $L_{d,k+1}\ge R_{d,k+1}$,
\item  if $\frac{d-2}{4}\le k\le \frac{d-2}{3}$, then $L_{d,k+1}\ge R_{d,k+1}$,
\item  if $\frac{d}{3}\le k\le \frac{2d-5}{3}$, then $L_{d,k+1}\le R_{d,k+1}$,
\item  if $\frac{2d-5}{3}\le k$, then $L_{d,k+1}\le R_{d,k+1}$,
\end{enumerate}
and the result follows. Indeed,
\begin{enumerate}
\item the factor in $L_{d,k+1}$  is positive and so $L_{d,k+1}\ge 2(d-3k-1)\ge R_{d,k+1}$,
\item the factor in $L_{d,k+1}$  is at least 1,  the exponent $d-k-1$ in
$L_{d,k+1}$ is at least  $\frac{d+2}{3}$  and $\frac{d+2}{2}\ge R_{d,k+1}$, so we have $L_{d,k+1}\ge 2^{\frac{d+2}{3}}(d-3k-1){\ge 2^{\frac{d+2}{3}}}  \ge\frac{d+2}{2}\ge R_{d,k+1}$,
\item the factor in $L_{d,k+1}$  is no more than  $-1$, the exponent $d-k-1$ in
$L_{d,k+1}$ is at least  $\frac{d+2}{3}$  and $ R_{d,k+1}\ge\frac{-d+10}{3}$, so we get
$R_{d,k+1}\ge\frac{-d+10}{3} \ge -2^{\frac{d+2}{3}}\ge2^{\frac{d+2}{3}}(d-3k-1)\ge L_{d,k+1}$,
\item  the factor in $L_{d,k+1}$  is non--positive and so $R_{d,k+1}\ge d-3k-1\ge L_{d,k+1}$.
\end{enumerate}

It follows that the change in the monotonicity of the sequence $I_{d,k}$ occurs in the interval  $\Z\cap [\frac{d-2}{3}, \frac{d}{3}]$.
For fixed $d\ge2$, we have found   the maximum in $I_{d,k}$ attained at $k=\lfloor\frac{d}{3}\rfloor=\begin{cases}
\frac{d}{3},& d\equiv0 \mod 3,\\
\frac{d-1}{3},& d\equiv1 \mod 3,\\
\frac{d-2}{3},& d\equiv2 \mod 3.
\end{cases}$
\endproof
\begin{cor}[B\'{a}r\'{a}ny conjecture holds for isocanted]
If $d\ge{2}$  and $0\le k<k+1\le d-1$, then
$I_{d,k}\ge\min\{I_{d,0}, I_{d,d-1}\}=I_{d,d-1}=(d+1)d$.
\end{cor}
\proof
Use unimodality and Proposition \ref{prop:extremes}.
\endproof

The $3^d$ conjecture and the flag conjecture were posed by Kalai in 1989, for centrally symmetric polytopes.

\begin{cor}[$3^d$ conjecture holds for isocanted]\label{cor:3d_conj}
For $d\ge2$, it holds $ \sum_{k=0}^{d}I_{d,k}=3^{d+1}-2^{d+2}+{2}$ and this is larger than $3^{d}$.
\end{cor}
\proof
  The binomial theorem $(x+y)^d= \sum_{j=0}^{d} x^j y^{d-j} {d \choose j}$ with $x=1$ yields $2^d=\sum_{j=0}^{d} {d \choose j}$ and $3^d=\sum_{j=0}^{d} 2^{d-j} {d \choose j}$. Then
\begin{equation}
3^{d+1}-2\times 2^{d+1}=\sum_{j=0}^{d+1} 2^{d+1-j} {{d+1} \choose j}-2 \sum_{j=0}^{d+1} {{d+1} \choose j}
=\sum_{j=0}^{d+1} (2^{d+1-j}-2) {{d+1} \choose j}=\sum_{j=0}^{d-1} I_{d,j} +\text{two summands.}
\end{equation}
Summand for  $j=d$ is zero and  for $j=d+1$ is ${-1}$, whence, using $I_{d,d}=1$, we get
the claimed equality.  Proof of the inequality:   we have $2^3=8=3^2-1$ and $2^{d-2}<3^{d-2}$. Multiply termwise and
get $2^{d+1}\le 3^{d-2}(3^2-1)=3^d-3^{d-2}<3^{d}+1$
whence $2(2^{d+1}-1)<2\times 3^d=3^{d+1}-3^d$.
\endproof

\begin{rem}
Recall that \emph{Stirling number of the second kind} is the number of ways to partition $[d]$ into $k$ non--empty subsets, and it is denoted $S(d,k)$.
We have $3^{d+1}-2^{d+2}+2=2S(d+2,3)+1$ (see Wikipedia and OEIS A101052, OEIS A028243 and OEIS A000392 in  \cite{OEIS}).
\end{rem}

\begin{rem}
Recall that a \emph{Hanner polytope} is obtained from closed intervals, by  using two operations any finite number of times:
Cartesian product and polar. They were studied by Hanner in 1956.
Is $\II^d$ a Hanner polytope? Conversely, is some  Hanner polytope an IAP?
Since Hanner polytopes satisfy the  $3^d$ conjecture and they attain the minimal conjectured value (see \cite{Sanyal_al}), then the answer is NO in both cases.
\end{rem}

Recall that  a \emph{complete flag} in a polytope $\PP$ is a maximal chain of faces of $\PP$ with increasing dimensions.
Next, we  count 
complete flags (and call them flags, for short). The number of  flags in a $d$--box is $2^dd!$, because there are $2^d$ vertices and,
at each one, there are $d!$  flags. The flag conjecture   yields
that boxes minimize flags among centrally symmetric polytopes; see \cite{Kalai,Sanyal_al,Senechal_Libro}.

\begin{cor}[Flag conjecture holds for isocanted]\label{cor:flag_conj}
The number of  flags in $\II^d$ is $(d+1)(d-1)!(2^{d+1}-4)$ and it is larger than $2^dd!$, for $d\ge2$.
\end{cor}

\proof
In $\II^d$  there are  $2(d+1)$ vertices of valence $d$, and the remaining $2(2^d-d-2)$ vertices have valence $d+1$. Indeed, the
vertices  of length 1 or $d$ have valence $d$. A vertex of length $2\le t\le d-1$ has valence $d+1$, because it has $t$ parents
and $d+1-t$ children. Reasoning as in boxes, we find  $d!$ flags
beginning at a vertex of
valence $d$. Using Item \ref{item:lattice} in Corollary \ref{cor:face_lattice;unique;cubical}, we find  $(d+1)(d-1)!$ flags beginning at
a vertex $V$ of valence $d+1$, because $\II^d$ is cubical and there are $d+1$ $(d-1)$--cuboids meeting at $V$.
Thus, adding up,
$2(d+1)\times d!+2(2^d-d-2)\times (d+1)(d-1)!=(d+1)(d-1)!(2^{d+1}-4)$ is the total number of flags.  Further, we have $(2^{d-1}-1)(d+1)>2^{d-2}d$, for $d\ge2$, whence the claimed inequality.
\endproof

The \emph{cubical lower bound conjecture  (CLBC)} was posed by Jockusch in 1993 and rephrased, in terms of the  \emph{cubical $g$--vector} $g^c$, by Adin et al. in 2019 as follows: is $g^c_{d,2}\ge0$?; see \cite{Adin_et_al,Jockusch}.

\begin{prop}[CLBC holds for isocanted]\label{prop:CLBC}
$g^c_{d,2}\ge0$ holds true for $\II^d$, for $d\ge2$.
\end{prop}
\proof
We have computed the sequence $g^c_{d,2}$ for IAPs,  obtaining   $6, 20, 50, 112, 238, \ldots$; see OEIS A052515 in \cite{OEIS}.
\endproof

Recall that the \emph{distance} between two vertices of a polytope is the minimum number of edges in a path joining them. The \emph{diameter}
of a polytope is the greatest distance between two vertices of the polytope.
\begin{cor}[Diameter of isocanted]\label{cor:diameter}
The diameter of $\II^d$ is $d+1$.
\end{cor}
\proof
Consider different proper subsets $W,W'\subset [d+1]$ and assume $|W\cap W'|=i$, $|W|=i+w$ ,$|W'|=i+w'$, with $i,w,w'\ge0$ and
$i+w+w'\le d+1$. To go from vertex $\underline{W}$ to vertex $\underline{W'}$ one must drop (one at a time) the $w$ digits in
$W\setminus W'$ and one must gain (one at a time)
the $w'$ digits in $W'\setminus W$, whence $\dd(\underline{W},\underline{W'})=w+w'$. In the particular case that $W'$ is complementary to $W$, we get the greatest distance $\dd(\underline{W},\underline{W'})=d+1$.
\endproof

\section{Future work}
We would like to  compute the
$f$--vector of  a general alcoved polytope.
\section*{Acknowledgments}\label{sec:acknowledgments} We thank the referee for a careful revision.

{\small
{\em Authors' addresses}:
{\em M.J.~de la Puente},  Universidad Complutense, Madrid,
Spain, e-mail:
 \texttt{mpuente@\allowbreak ucm.es}.
 {\em P.L.~Claver\'{\i}a},  Universidad de Zaragoza, Zaragoza,
Spain, e-mail:
 \texttt{plcv@\allowbreak unizar.es}.
 }

\end{document}